\documentclass[10pt,draft,twoside,leqno,a4paper]{article}
\usepackage{amsmath,amsthm,amsfonts,amssymb,verbatim}
\newtheorem{Theorem}{Theorem}
\newtheorem{Lemma}{Lemma}

\newtheorem{Proposition}{Proposition}

\newcommand{\dash}{\mathchoice
    {\mkern.70mu\raise.50ex\hbox{\vrule height.1ex width.40em depth0pt}}
    {\mkern.40mu\raise.48ex\hbox{\vrule height.1ex width.30em depth0pt}}
    {\mkern.33mu\raise.30ex\hbox{\vrule height.1ex width.25em depth0pt}}
    {\mkern.10mu\raise.20ex\hbox{\vrule height.1ex width.20em depth0pt}}
    }

\def\esssup{\mathop{\mathrm {ess\,sup}}}

\newcommand{\R}{\mathbb R}
\newcommand{\C}{\mathbb C}

\newcommand{\de}{\partial}
\newcommand{\dbar}{\overline\partial}

\newcommand{\bra}{\langle}
\newcommand{\ket}{\rangle}

\newcommand{\fo}{f_0}

\newcommand{\ko}{k_0}
\newcommand{\mo}{\mu_0}

\renewcommand{\d}{{\mathrm d}}

\renewcommand{\Re}{\mathrm{Re}}
\renewcommand{\Im}{\mathrm{Im}}
\begin{document}
\title{On Beltrami equations and H\"older regularity 
}
\author{Tonia Ricciardi\thanks{
Supported in part by 
the MIUR National Project {\em Variational Methods and
Nonlinear Differential Equations}.}\\
{\small Dipartimento di Matematica e Applicazioni
``R.~Caccioppoli"}\\
{\small Universit\`a di Napoli Federico II}\\
{\small Via Cintia, 80126 Napoli, Italy}\\
{\small fax: +39 081 675665}\\
{\small\tt{ tonia.ricciardi@unina.it}}\\
}
\date{}
\maketitle
\begin{abstract}
We estimate the H\"older exponent $\alpha$ of solutions to 
the Beltrami equation $\dbar f=\mu\de f$,
where the Beltrami coefficient satisfies $\|\mu\|_\infty<1$. 
Our estimate improves the classical estimate $\alpha\ge\|K_\mu\|^{-1}$,
where $K_\mu=(1+|\mu|)/(1-|\mu|)$, and it 
is sharp, in the sense that it is actually attained 
in a class of mappings which generalize the radial stretchings.
Some other properties of such mappings are also provided.
\end{abstract}
\begin{description}
\item {\textsc{Key Words:}} linear Beltrami equation, measurable coefficients,
H\"older regularity
\item {\textsc{MSC 2000 Subject Classification:}} 30C62 (35J25)
\end{description}
\section{Introduction and main results}
Let $\Omega$ be a bounded open subset of $\R^2$ and let
$f\in W_{\mathrm{loc}}^{1,2}(\Omega,\C)$ satisfy the Beltrami equation
\begin{equation}
\label{beltrami}
\dbar f=\mu\,\de f
\qquad\mathrm{a.e.\ in\ }\Omega,
\end{equation}
where $\mu\in L^\infty(\Omega,\C)$ satisfies $\|\mu\|_\infty<1$,
$\dbar=(\de_1+i\de_2)/2$, $\de=(\de_1-i\de_2)/2$.
By classical results, see, e.g.,
Iwaniec and Martin~\cite{IM}, 
there exists $\alpha\in(0,1)$ such that $f$ is $\alpha$-H\"older continuous 
in $\Omega$. Namely, for every compact $T\Subset\Omega$ there exists $C_T>0$
such that
\[
\frac{|f(z)-f(z')|}{|z-z'|^\alpha}\le C_T
\qquad\forall z,z'\in T,\ z\neq z'.
\]
Moreover, let 
\[
K_\mu=\frac{1+|\mu|}{1-|\mu|} 
\]
denote the distortion function. 
Then, the following estimate holds:
\begin{equation}
\label{classicalhoelder}
\alpha\ge\frac{1}{\|K_\mu\|_\infty}.
\end{equation}
This estimate is sharp, in the sense that it reduces to
an equality on the radial stretchings
\begin{equation}
\label{radialstretch}
f(z)=|z|^{\alpha-1}z.
\end{equation}
Indeed, the function $f$ defined above satisfies the Beltrami equation~\eqref{beltrami}
with Beltrami coefficient given by $\mu(z)=-(1-\alpha)/(1+\alpha)z\bar z^{-1}$, and a simple
computation yields $K_\mu\equiv\alpha^{-1}$.
We also note that $f$ defined in \eqref{radialstretch}
satisfies the pointwise equality
\begin{equation}
\label{qc}
|Df(z)|^2=\alpha^{-1}J_f,
\end{equation}
which is of interest in the context of quasiconformal mappings.
\par
There exists a wide literature concerning the regularity theory
for \eqref{beltrami}, in the degenerate case where
$\|\mu\|_\infty=1$, or equivalently, when the distortion function $K_\mu$ is unbounded. 
See, e.g., \cite{AIKM,GMSV,IM,IS}, and the references therein.
Our aim in this note is to extend the classical estimate \eqref{classicalhoelder}
in another sense, by subtracting a correction term to $K_\mu$.
\begin{Theorem}
\label{thm:estimate}
Let $f\in W_{\mathrm{loc}}^{1,2}(\Omega,\C)$ satisfy the Beltrami equation~\eqref{beltrami}.
Then, $f$ is $\alpha$-H\"older continuous with
\begin{equation}
\label{alphaestimate}
\alpha\ge\left(\esssup_{S_\rho(x)\subset\Omega}
\frac{1}{|S_\rho(x)|}\int_{S_\rho(x)}\Big\{K_\mu
-2\frac{|\mu|+\Re\left(\mu,n^2\right)}{1-|\mu|^2}\Big\}\right)^{-1},
\end{equation}
where $n$ denotes the outer unit normal. 
\end{Theorem}
Here, $(\cdot,\cdot)$ denotes the standard inner product in $\C$,
and $n^2$ is understood in the complex multiplication sense.
That is, for every $y=x+\rho e^{it}\in S_\rho(x)$, we have
$n(y)=e^{it}$, $(\mu,n^2)(y)=\mu(y)\overline{e^{2it}}=\mu(y) e^{-2it}$.
In particular, $|\mu|+\Re\left(\mu,n^2\right)\ge0$, and therefore
\begin{align*}
\esssup_{S_\rho(x)\subset\Omega}
\frac{1}{|S_\rho(x)|}&\int_{S_\rho(x)}\Big\{K_\mu
-2\frac{|\mu|+\Re\left(\mu,n^2\right)}{1-|\mu|^2}\Big\}\\
\le&\esssup_{S_\rho(x)\subset\Omega}
\frac{1}{|S_\rho(x)|}\int_{S_\rho(x)}K_\mu\le\|K_\mu\|_\infty.
\end{align*}
It follows that estimate \eqref{alphaestimate} indeed improves \eqref{classicalhoelder}.
Moreover, our estimates are sharp, in the following sense.
For every $M>1$, let 
\begin{align}
\label{c}
c=c(M)=\frac{2}{1+M^{-1}}
\end{align}
and let $\ko:\R\to\R$ be the bounded, $2\pi$-periodic function defined by
\begin{equation}
\label{ko}
\ko(\theta)=\begin{cases}
1,&\mathrm{if\ }\theta\in[0,c\pi/2)\cup[\pi,\pi+c\pi/2)\\
M,&\mathrm{otherwise}.
\end{cases}
\end{equation}
Let $\Theta_1,\Theta_2:\R\to\R$ be the $2\pi$-periodic Lipschitz functions defined by
\[
\Theta_1(\theta)=\begin{cases}
\sin(c^{-1}\theta-\pi/4),&\theta\in[0,c\pi/2)\\
\cos(c^{-1}M(\theta-c\pi/2)-\pi/4),&\theta\in[c\pi/2,\pi)\\
-\sin(c^{-1}(\theta-\pi)-\pi/4),&\theta\in[\pi,\pi+c\pi/2)\\
-\cos(c^{-1}M(\theta-\pi-c\pi/2)-\pi/4),&\theta\in[\pi+c\pi/2,2\pi)
\end{cases}
\]
and
\[
\Theta_2(\theta)=\begin{cases}
-\cos(c^{-1}\theta-\pi/4),&\theta\in[0,c\pi/2)\\
\sin(c^{-1}M(\theta-c\pi/2)-\pi/4),&\theta\in[c\pi/2,\pi)\\
\cos(c^{-1}(\theta-\pi)-\pi/4),&\theta\in[\pi,\pi+c\pi/2)\\
-\sin(c^{-1}M(\theta-\pi-c\pi/2)-\pi/4),&\theta\in[\pi+c\pi/2,2\pi).
\end{cases}
\]
The sharpness of Theorem~\ref{thm:estimate} is a consequence of the following.
\begin{Theorem}
\label{thm:sharp}
Let $B$ the unit disk in $\R^2$ and let $\fo\in W^{1,2}(B,\C)$ be defined 
in $B\setminus\{0\}$ by
\[
\fo(z)=|z|^{1/c}\left(\Theta_1(\arg z)+i\Theta_2(\arg z)\right).
\]
Let
\[
\mo(z)=\frac{1-\ko(\arg z)}{1+\ko(\arg z)}z\bar z^{-1}.
\]
Then $\fo$ satisfies \eqref{beltrami} with $\mu=\mo$.
Furthermore, there exists $\bar M>1$ such that
\begin{align*}
\esssup_{S_\rho(x)\subset\Omega}
\frac{1}{|S_\rho(x)|}&\int_{S_\rho(x)}\Big\{K_{\mo}
-2\frac{|\mo|+\Re\left(\mo,n^2\right)}{1-|\mo|^2}\Big\}=c,
\end{align*}
for every $M\in(1,\bar M)$.
\end{Theorem}
Our proof of Theorem~\ref{thm:estimate} is based on the equivalence
between Beltrami equations and elliptic divergence form equations 
with unit determinant, to which we can apply some recent results in \cite{besthoelder}.
On the other hand, in order to prove Theorem~\ref{thm:sharp}, we shall use some properties of
a class of mappings which generalize the radial stretchings.
More precisely, we consider mappings of the form:
\begin{equation}
\label{fansatz}
f(z)=|z|^\alpha(\eta_1(\arg z)+i\eta_2(\arg z))
\end{equation}
and coefficients of the form
\begin{equation}
\label{mu}
\mu(z)=\frac{1-k(\arg z)}{1+k(\arg z)}z\bar z^{-1}
\end{equation}
for some $2\pi$-periodic functions 
$\eta_1,\eta_2,k$,
$k\ge1$.
We prove:
\begin{Proposition}
\label{prop:polar}
Suppose $\mu$ is of the form \eqref{mu} and $f$ satisfies \eqref{fansatz}.
Then $f$ is a a solution to \eqref{beltrami} if and only if $(\eta_1,\eta_2)$ 
satisfies the first order system:
\begin{align}
\label{system}
\begin{cases}
\eta_1'=-\alpha k\eta_2\\
\eta_2'=\alpha k\eta_1.
\end{cases}
\end{align}
Furthermore, for every $z\neq0$ 
the following equality holds:
\begin{equation}
\label{quasi}
|Df|^2=k(\arg z)J_f,
\end{equation}
where $|Df|$ denotes the operator norm of $Df$.
\end{Proposition}
We note that mappings of the form \eqref{fansatz} generalize the radial stretchings
\eqref{radialstretch}.
Indeed, when $\eta_1(\theta)=\cos\theta$ and $\eta_2(\theta)=\sin\theta$,
\eqref{fansatz} reduces to \eqref{radialstretch}.
We expect that mappings of the form \eqref{fansatz} should be of interest 
in relation to other results on quasiconformal mapping theory and elliptic equations,
as well.
\subparagraph*{Notation}
Henceforth, for every measurable function $f$, $\sup f$ denotes the essential upper bound
of $f$. All integrals are taken with respect to the Lebesgue measure.
\section{Proofs}
\label{sec:proofs}
We use some results in
\cite{besthoelder} for solutions to the elliptic divergence form equation
\begin{equation}
\label{elliptic}
\mathrm{div}(A\nabla\cdot)=0\qquad\mathrm{in\ }\Omega
\end{equation}
where $A$ is a bounded, symmetric matrix-valued function satisfying
\begin{equation}
\label{det1}
\det A=1\qquad\mathrm{\ a.e.\ in\ }\Omega. 
\end{equation}
More precisely, let 
\begin{equation}
\label{J}
J(\theta)=\left(\begin{matrix}
&\cos\theta&&-\sin\theta\\
&\sin\theta&&\cos\theta
\end{matrix}\right).
\end{equation}
The following results were established in \cite{besthoelder}.
\begin{Theorem}[\cite{besthoelder}]
\label{thm:besthoelder}
The following estimates hold.
\begin{enumerate}
\item[(i)]
Let $w\in W_{\mathrm{loc}}^{1,2}(\Omega)$ be a weak solution to \eqref{elliptic},
where $A$ satisfies \eqref{det1}.
Then, $w$ is $\alpha$-H\"older continuous with
\[
\alpha\ge\left(\sup_{S_\rho(x)\subset\Omega}
\frac{1}{|S_\rho(x)|}\int_{S_\rho(x)}
\bra n,A n\ket\right)^{-1},
\]
where $n$ denotes the outer unit normal.
\item[(ii)]
Let $A_0$ be the symmetric matrix-valued function satisfying \eqref{det1}
defined for every $z\neq0$ by
\[
A_0(z)=J(\arg z)
\left(\begin{matrix}
&\ko(\arg z)&&0\\&0&&\ko^{-1}(\arg z)
\end{matrix}\right)
J^*(\arg z),
\]
where $\ko$ is the function defined in \eqref{ko}.
There exists $\bar M>1$ such that
\begin{equation}
\label{supso}
\sup_{S_\rho(x)\subset\Omega}
\frac{1}{|S_\rho(x)|}\int_{S_\rho(x)}\bra n,A_0 n\ket
=c=\frac{2}{1+M^{-1}},
\end{equation}
for every $M\in(1,\bar M)$.
Furthermore, the function $u_0=|z|^{1/c}\Theta_1(\arg z)$ 
is a weak solution to \eqref{elliptic} with $A=A_0$.
\end{enumerate}
\end{Theorem}
The reduction of 
Beltrami equations to
elliptic divergence form equations satisfying \eqref{det1}
is well-known.
For the reader's convenience and for uniformity of notation, we sketch it below.
\begin{Lemma}
\label{lem:beltramitoelliptic}
Let $f\in W_{\mathrm{loc}}^{1,2}(\Omega,\C)$ satisfy the Beltrami equation~\eqref{beltrami}.
Let
$A_\mu=(a_{\mu,ij})$ be the bounded, symmetric matrix-valued function
satisfying \eqref{det1}, 
defined in terms of the Beltrami coefficient $\mu=\mu_1+i\mu_2$, $\mu_1(z),\mu_2(z)\in\R$, 
by
\begin{align}
\label{Amu}
&a_{\mu,11}=\frac{1-2\mu_1+|\mu|^2}{1-|\mu|^2}
&&a_{\mu,22}=\frac{1+2\mu_1+|\mu|^2}{1-|\mu|^2}\\
\nonumber
&a_{\mu,12}=a_{\mu,21}=-\frac{2\mu_2}{1-|\mu|^2}.
\end{align}
Then $u=\Re f$ and $v=\Im f$ are a weak solutions to the elliptic equation
\eqref{elliptic}
with $A=A_\mu$.
\end{Lemma}
\begin{proof}
We have:
\begin{align*}
&\dbar f=\frac{1}{2}\left(\begin{matrix}u_x-v_y\\u_y+v_x\end{matrix}\right)
&&\de f=\frac{1}{2}\left(\begin{matrix}u_x+v_y\\-u_y+v_x\end{matrix}\right).
\end{align*}
Setting
\begin{align*}
&Q=\left(\begin{matrix}&0&&-1\\&1&&0\end{matrix}\right)
&&R=\left(\begin{matrix}&1&&0\\&0&&-1\end{matrix}\right),
\end{align*}
for every $z=a+ib=(a,b)^T$ we have
\begin{align*}
&Qz=\left(\begin{matrix}-b\\a\end{matrix}\right)=iz,
&&Rz=\left(\begin{matrix}a\\-b\end{matrix}\right)=\overline z.
\end{align*}
Hence, we can write
\begin{align*}
&\dbar f=\frac{1}{2}\left(\nabla u+Q\nabla v\right),
&&\de f=\frac{R}{2}\left(\nabla u-Q\nabla v\right).
\end{align*}
Setting
\[
M=\left(\begin{matrix}&\mu_1&&-\mu_2\\&\mu_2&&\mu_1\end{matrix}\right),
\]
equation~\eqref{beltrami} can be written in the form:
\begin{equation*}
\nabla u+Q\nabla v=MR\left(\nabla u-Q\nabla v\right).
\end{equation*}
It follows that
\[
\left(I-MR\right)\nabla u=-\left(I+MR\right)Q\nabla v
\]
and consequently $u$ satisfies
\[
\left(I+MR\right)^{-1}\left(I-MR\right)\nabla u=-Q\nabla v
\]
and $v$ satisfies
\[
-Q\left(I-MR\right)^{-1}\left(I+MR\right)Q\nabla v=Q\nabla u.
\]
A straightforward computation yields
\[
\left(I+MR\right)^{-1}\left(I-MR\right)
=A_\mu=-Q\left(I-MR\right)^{-1}\left(I+MR\right)Q.
\]
Now the conclusion follows observing that $\mathrm{div}(Q\nabla\cdot)=0$.
\end{proof}
\begin{proof}[Proof of Theorem~\ref{thm:estimate}]
In view of Lemma~\ref{lem:beltramitoelliptic} and Theorem~\ref{thm:besthoelder},
$\Re f$ and $\Im f$ are $\alpha$-H\"older continuous with
\[
\alpha\ge\left(\sup_{S_\rho(x)\subset\Omega}
\frac{1}{|S_\rho(x)|}\int_{S_\rho(x)}
\bra n,A_\mu n\ket\right)^{-1},
\]
where $A_\mu$ is defined in \eqref{Amu}.
Setting $\xi=x+\rho e^{it}$, $t\in\R$ for every $\xi\in S_\rho(x)\subset\Omega$, 
we have $n(\xi)=e^{it}$.
Hence, we compute
\begin{align*}
\bra n(\xi),A_\mu(\xi)n(\xi)\ket=\bra e^{it},A_\mu(\xi)e^{it}\ket
=&a_{\mu,11}\cos^2t+2a_{\mu,12}\sin t\cos t+a_{\mu,22}\sin^2t\\
=&\frac{1+|\mu|^2}{1-|\mu|^2}-2\frac{\mu_1\cos2t+\mu_2\sin2t}{1-|\mu|^2}.
\end{align*}
Using the identity
\begin{align*}
\frac{1+|\mu|^2}{1-|\mu|^2}+2\frac{|\mu|}{1-|\mu|^2}=\frac{1+|\mu|}{1-|\mu|}=K_\mu,
\end{align*}
we have on $S_\rho(x)$
\begin{equation}
\label{AK}
\bra n,A_\mu n\ket=K_\mu
-2\frac{|\mu|+\Re\left(\mu,n^2\right)}{1-|\mu|^2}
\end{equation}
and we conclude the proof.
\end{proof}
Now we prove some properties for the special case
where the Beltrami coefficient $\mu$ is of the form~\eqref{mu}.
\begin{Lemma}
\label{lem:mutoA}
Suppose $\mu$ is of the form \eqref{mu}.
Then $A_\mu$ as defined in \eqref{Amu} is given by
\begin{align*}
A_\mu(z)=&J\left(
\begin{matrix}&k(\arg z)&&0\\
&0&&k^{-1}(\arg z)
\end{matrix}\right)J^*\\
=&\left(
\begin{matrix}&k\cos^2\theta+k^{-1}\sin^2\theta &&\left(k-k^{-1}\right)\sin\theta\cos\theta\\
&\left(k-k^{-1}\right)\sin\theta\cos\theta&&k^{-1}\cos^2\theta+k\sin^2\theta 
\end{matrix}\right).
\end{align*} 
\end{Lemma}
\begin{proof}
Setting $\theta=\arg z$, we have
\[
\mu=\frac{1-k}{1+k}e^{2i\theta}
\]
and therefore,
\begin{align*}
&\mu_1=\frac{1-k}{1+k}\cos2\theta&&\mu_2=\frac{1-k}{1+k}\sin2\theta.
\end{align*}
We compute
\begin{align*}
&1-|\mu|^2=\frac{4k}{(1+k)^2}
&&1+|\mu|^2=2\frac{1+k^2}{(1+k)^2}.
\end{align*}
In view of \eqref{Amu}, we have
\begin{align*}
(1-|\mu|^2)a_{\mu,11}=&1-2\mu_1+|\mu|^2=2\frac{1-\cos2\theta+k^2(1+\cos2\theta)}{(1+k)^2}\\
=&4\frac{\sin^2\theta+k^2\cos^2\theta}{(1+k)^2}
\end{align*}
and consequently
\[
a_{\mu,11}=k\cos^2\theta+\frac{1}{k}\sin^2\theta.
\]
Similarly,
\[
a_{\mu,22}=\frac{1}{k}\cos^2\theta+k\sin^2\theta.
\]
Finally, 
\begin{align*}
a_{\mu,12}=&-2\frac{\mu_2}{1-|\mu|^2}\\
=&-2\frac{1-k}{1+k}\sin2\theta\frac{(1+k)^2}{4}=\left(k-\frac{1}{k}\right)\sin\theta\cos\theta.
\end{align*}
\end{proof}
In what follows it will be convenient to use polar cooordinates
$x=\rho\cos\theta$, $y=\rho\sin\theta$. We recall that 
\begin{align*}
\dbar=&\frac{1}{2}(\de_x+i\de_y)=\frac{1}{2}\left(\frac{z}{|z|}\de_\rho
+i\frac{z}{|z|^2}\de_\theta\right)\\
\de=&\frac{1}{2}(\de_x-i\de_y)=\frac{1}{2}\left(\frac{\bar z}{|z|}\de_\rho
-i\frac{\bar z}{|z|^2}\de_\theta\right).
\end{align*}
Hence, \eqref{beltrami} is equivalent to
\[
\frac{z}{|z|}f_\rho+i\frac{z}{|z|^2}f_\theta
=\mu\left(\frac{\bar z}{|z|}f_\rho-i\frac{\bar z}{|z|^2}f_\theta\right),
\]
from which we derive the polar form of the Beltrami equation~\eqref{beltrami}:
\begin{equation}
\label{fpolar}
(e^{i\theta}-\mu e^{-i\theta})f_\rho=-\frac{i}{\rho}(e^{i\theta}+\mu e^{-i\theta})f_\theta.
\end{equation}
Finally, we prove Proposition~\ref{prop:polar}.
\begin{proof}[Proof of Proposition~\ref{prop:polar}]
Proof of \eqref{system}.
If $\mu$ and $f$ are as given, we have
\begin{align*}
e^{i\theta}-\mu e^{-i\theta}
=&e^{i\theta}\left(1-\frac{1-k}{1+k}\right)=\frac{2k}{1+k}e^{i\theta}\\
e^{i\theta}+\mu e^{-i\theta}
=&e^{i\theta}\left(1+\frac{1-k}{1+k}\right)=\frac{2}{1+k}e^{i\theta}.
\end{align*} 
Since $f_\rho=\alpha\rho^{\alpha-1}(\eta_1+i\eta_2)$
and $f_\theta=\rho^\alpha(\eta_1'+i\eta_2')$,
we derive from \eqref{fpolar}
\[
\alpha k(\eta_1+i\eta_2)
=-i(\eta_1'+i\eta_2'),
\]
which in turn implies \eqref{system}.
\par
Proof of \eqref{quasi}.
We regard $f$ as a mapping from $(\Omega,\d\rho^2+\rho^2\d\theta^2)$
to $(\R^2,\d x^2+\d y^2)$.
Writing  $f=\rho^\alpha\eta_1\de_1+\rho^\alpha\eta_2\de_2$, 
for every $\xi=\xi_\rho\de_\rho+\xi_\theta\de_\theta$
we have
\[
\d f\xi=\rho^{\alpha-1}\left[(\alpha\eta_1\xi_\rho+\rho\eta_1'\xi_\theta)\de_1
+(\alpha\eta_2\xi_\rho+\rho\eta_2'\xi_\theta)\de_2\right].
\]
Hence,
\[
\rho^{-2(\alpha-1)}|\d f\xi|^2=\alpha^2(\eta_1^2+\eta_2^2)\xi_\rho^2
+2\alpha\rho(\eta_1\eta_1'+\eta_2\eta_2')\xi_\rho\xi_\theta
+\rho^2(\eta_1'^2+\eta_2'^2)\xi_\theta^2.
\]
In view of \eqref{system}, we have
$\eta_1\eta_1'+\eta_2\eta_2'=0$ and $\eta_1'^2+\eta_2'^2=\alpha^2k^2(\eta_1^2+\eta_2^2)$.
Therefore,
\[
|\d f\xi|^2=\rho^{2(\alpha-1)}\alpha^2(\eta_1^2+\eta_2^2)
\left(\xi_\rho^2+\rho^2k^2\xi_\theta^2\right).
\]
On the other hand, $|\xi|^2=\xi_\rho^2+\rho^2\xi_\theta^2$.
It follows that
\begin{align*}
|\d f|^2=&\sup_{\xi\neq0}
\frac{\alpha^2\rho^{2(\alpha-1)}(\eta_1^2+\eta_2^2)(\xi_\rho^2+\rho^2k^2\xi_\theta^2)}
{\xi_\rho^2+\rho^2\xi_\theta^2}\\
=&\alpha^2\rho^{2(\alpha-1)}(\eta_1^2+\eta_2^2)\sup_{t>0}\frac{1+\rho^2k^2t}{1+\rho^2t}
=\alpha^2k^2\rho^{2(\alpha-1)}(\eta_1^2+\eta_2^2).
\end{align*}
In order to compute $J_f$, we note that in polar coordinates on $\R^2$
the Hodge star  operator satisfies $\ast1=\rho\d\rho\wedge\d\theta$ and $\ast\ast=1$.
Hence,
\begin{align*}
J_f=\ast(\d\Re f\wedge\d\Im f)
=&\alpha\rho^{2\alpha-1}(\eta_1\eta_2'-\eta_1'\eta_2)\ast(\d\rho\wedge\d\theta)\\
=&\alpha\rho^{2(\alpha-1)}(\eta_1\eta_2'-\eta_1'\eta_2).
\end{align*}
In view of \eqref{system}, we have $\eta_1\eta_2'-\eta_1'\eta_2=\alpha k(\eta_1^2+\eta_2^2)$.
It follows that
\[
J_f=\alpha^2k\rho^{2(\alpha-1)}(\eta_1^2+\eta_2^2)
\]
and, finally, that
\[
|\d f|^2=kJ_f
\]
for every $z\neq0$. Now \eqref{quasi} is established.
\end{proof}
\begin{proof}[Proof of Theorem~\ref{thm:sharp}]
By direct check, $(\Theta_1,\Theta_2)$ is a solution to 
system~\eqref{system} with $k=\ko$ and $\alpha=c^{-1}$,
where $c$ is defined in \eqref{c}.
Hence, in view of Proposition~\ref{prop:polar},
the function $f_0=\rho^{1/c}(\Theta_1+i\Theta_2)$ is a solution
to the Beltrami equation~\eqref{beltrami} with $\mu=\mo$.
By Lemma~\ref{lem:beltramitoelliptic}, $\Re\fo$ and $\Im\fo$
are solutions to the elliptic equation~\eqref{elliptic},
with $A=A_0$ the matrix-valued funtion defined in Theorem~\ref{thm:besthoelder}--(ii).
In view of \eqref{AK}, we have
\begin{align*}
c=&\left(\sup_{S_\rho(x)\subset\Omega}\int_{S_\rho(x)}\bra n,A_0 n\ket\right)^{-1}\\
=&\left(\sup_{S_\rho(x)\subset\Omega}\int_{S_\rho(x)}\Big\{K_{\mo}
-2\frac{|\mo|+\Re\left(\mo,n^2\right)}{1-|\mo|^2}\Big\}\right)^{-1}
\end{align*}
for every $M\in(1,\bar M)$.
\end{proof}
\section*{Acknowledgments}
I thank Professor Kari Astala for
interesting and fruitful discussions.

\end{document}